\input amstex
\documentstyle{amsppt}
\magnification=\magstep1                        
\hsize6.5truein\vsize8.9truein                  
\NoRunningHeads
\loadeusm

\magnification=\magstep1                        
\hsize6.5truein\vsize8.9truein                  
\NoRunningHeads
\loadeusm

\document
\topmatter

\title
Do flat skew-reciprocal Littlewood polynomials exist?  
\endtitle

\rightheadtext{}

\author Tam\'as Erd\'elyi
\endauthor

\address Department of Mathematics, Texas A\&M University,
College Station, Texas 77843, College Station, Texas 77843 \endaddress

\thanks {{\it 2010 Mathematics Subject Classifications.} 11C08, 41A17, 26C10, 30C15}
\endthanks

\keywords
\endkeywords

\date July 6, 2020
\enddate

\email terdelyi\@math.tamu.edu
\endemail

\abstract
Polynomials with coefficients in $\{-1,1\}$ are called Littlewood polynomials.
Using special properties of the Rudin-Shapiro polynomials and classical results in approximation theory 
such as Jackson's Theorem, de la Vall\'ee Poussin sums, Bernstein's inequality, Riesz's Lemma, divided differences, 
etc., we give a significantly simplified  proof of a recent breakthrough result by Balister, Bollob\'as, Morris, 
Sahasrabudhe, and Tiba stating that there exist absolute constants $\eta_2 > \eta_1 > 0$ and a 
sequence $(P_n)$ of Littlewood polynomials $P_n$ of degree $n$ such that 
$$\eta_1 \sqrt{n} \leq |P_n(z)| \leq \eta_2 \sqrt{n}\,, \qquad z \in {\Bbb C}, \enskip |z| = 1\,,$$
confirming a conjecture of Littlewood from 1966. Moreover, the existence of a sequence $(P_n)$ of Littlewood
polynomials $P_n$ is shown in a way that in addition to the above flatness properties a certain symmetry is 
satisfied by the coefficients of $P_n$ making the Littlewood polynomials $P_n$ close to skew-reciprocal. 
\endabstract

\endtopmatter

\head 1. The Theorem \endhead

Polynomials with coefficients in $\{-1,1\}$ are called Littlewood polynomials.

\proclaim{Theorem 1.1}
There exist absolute constants $\eta_2 > \eta_1 > 0$ and a sequence $(P_n)$ of Littlewood polynomials $P_n$ 
of degree $n$ such that 
$$\eta_1 \sqrt{n} \leq |P_n(z)| \leq \eta_2 \sqrt{n}\,, \qquad z \in {\Bbb C}, \enskip |z| = 1\,. \tag 1.1$$
\endproclaim

Note that Beck [B-91] showed the existence of flat unimodular polynomials $P_n$ satisfying (1.1) with coefficients 
in the set of $k$th roots of unity. Beck showed the existence of flat unimodular polynomials $P_n$ 
of degree $n$ satisfying (1.1) with coefficients in the set of $k$th roots of unity and gave the value 
$k = 400$, but correcting a minor error in Beck's paper Belshaw [B-13] showed that the value of $k$ in [4] 
should have been $851$. Repeating Spencer's calculation Belshaw improved the value $851$ to $492$ in Beck's 
result, and an improvement of Spencer's method, due to Kai-Uwe Schmidt, allowed him to lower the value of
$k$ to $345$. The recent breakthrough result by Balister, Bollob\'as, Morris, Sahasrabudhe, and Tiba [B-20] 
formulated in Theorem 1.1 confirms a conjecture of Littlewood from 1966. Using special properties of the Rudin-Shapiro 
polynomials and classical results in approximation theory such as Jackson's Theorem, de la Vall\'ee Poussin sums, 
Bernstein's inequality, Riesz's Lemma, divided differences, etc., in this paper we give a significantly 
simplified proof of this beautiful and deep theorem. Moreover, the existence of a sequence $(P_n)$ of Littlewood 
polynomials $P_n$ is shown so that in addition to (1.1) a certain symmetry is satisfied by the coefficients of $P_n$.  

\proclaim{Theorem 1.2}
There exist absolute constants $0 < \eta_1 < \eta_2$, $\eta > 0$, and a sequence $(P_{2n})$ of Littlewood polynomials $P_{2n}$ 
of the form
$$P_{2n}(z) = \sum_{j=0}^{2n}{a_{j,n}z^j}\,, \qquad a_{j,n} \in \{-1,1\}\,, \enskip j=0,1,\ldots,2n, \enskip n=1,2,\ldots\,,$$ 
such that in addition to (1.1) the coefficients of $P_{2n}$ satisfy
$$a_{j,n} = -a_{2n-j,n}\,, \qquad 0 \leq j < n-m_n\,,$$
and
$$a_{j,n} = (-1)^{n-j}a_{2n-j,n}\,, \qquad n-m_n \leq j \leq n\,,$$
with some integers $0 \leq \eta n \leq m_n \leq n$.
\endproclaim

The theorem above may be viewed as a result in an effort to answer the following question.

\proclaim{Problem 1.3} 
Are there absolute constants $0 < \eta_1 < \eta_2$ and a sequence $(P_{2n})$ of skew-reciprocal 
Littlewood polynomials $P_{4n}$of the form
$$P_{4n}(z) = \sum_{j=0}^{2n}{a_{j,n}z^j}\,, \qquad a_{j,n} \in \{-1,1\}\,, \enskip j=0,1,\ldots,2n, \enskip n=1,2,\ldots\,,$$
such that in addition to (1.1) the coefficients of $P_{2n}$ satisfy
$$a_{j,n} = (-1)^{-j}a_{4n-j,n}\,, \qquad j=0,1,\ldots,2n\,?$$  
\endproclaim

This problem remains open. We remark that it is easy to see that every self-reciprocal Littlewood polynomial 
of the form 
$$P_{n}(z) = \sum_{j=0}^{n}{a_{j,n}z^j}\,, \qquad a_{j,n} \in \{-1,1\}\,, \enskip j=0,1,\ldots,n\,,$$
satisfying 
$$a_{j,n} = a_{n-j,n}\,, \qquad j=0,1,\ldots,n\,,$$ 
has at least one zero on the unit circle, see Theorem 2.8 in [E-11], or Corollary 2.5 in [M-06], for example. 
Hence there are no absolute constant $\eta_1 > 0$ and a sequence $(P_n)$ of self-reciprocal Littlewood polynomials $P_n$ 
of degree $n$ such that 
$$\eta_1 \sqrt{n} \leq |P_n(z)|\,, \qquad z \in {\Bbb C}, \enskip |z| = 1\,, \enskip n=1,2,\ldots\,.$$

\head 2. Rudin-Shapiro polynomials \endhead

Section 4 of [B-02] is devoted to the study of Rudin-Shapiro polynomials.
A sequence of Littlewood polynomials that satisfy just the upper bound of Theorem 1.1 is given 
by the Rudin-Shapiro polynomials. The Rudin-Shapiro polynomials appear in Harold Shapiro's 1951 
thesis [S-51] at MIT and are sometimes called just Shapiro polynomials. They also arise independently 
in Golay's paper [G-51]. 
The Rudin-Shapiro polynomials are remarkably simple to construct. They are defined recursively as follows:
$$\split P_0(z) & :=1\,, \qquad Q_0(z) := 1\,, \cr 
P_{m+1}(z) & := P_m(z) + z^{2^m}Q_m(z)\,, \cr
Q_{m+1}(z) & := P_m(z) - z^{2^m}Q_m(z)\,, \cr \endsplit$$
for $m=0,1,2,\ldots\,.$ Note that both $P_m$ and $Q_m$ are polynomials of degree $M-1$ with $M := 2^m$
having each of their coefficients in $\{-1,1\}$.
It is well known and easy to check by using the parallelogram law that
$$|P_{m+1}(e^{it})|^2 + |Q_{m+1}(e^{it})|^2 = 2(|P_m(e^{it})|^2 + |Q_m(e^{it})|^2)\,, \qquad t \in {\Bbb R}\,.$$
Hence
$$|P_m(e^{it})|^2 + |Q_m(e^{it})|^2 = 2^{m+1} = 2M\,, \qquad t \in {\Bbb R}\,. \tag 2.1$$
Observing that the first $2^m$ terms of $P_{m+1}$ are the same as the $2^m$ terms of $P_m$, we can define 
the polynomial $P_{<n}$ of degree $n-1$ so that its terms are the first $n$ terms of all $P_m$ for all $m$ for which 
$2^m \geq n$. The following bound, which is a straightforward consequence of (2.1) was proved by Shapiro [S-51]. 

\proclaim{Lemma 2.1}
We have
$$|P_{<n}(e^{it})| \leq 5\sqrt{n}\,, \qquad t \in {\Bbb R}\,.$$
\endproclaim

It is also well-known that 
$$P_m(1) = \|P_m(e^{it})\| := \max_{t \in {\Bbb R}}{|P_m(e^{it})}|= 2^{(m+1)/2}$$ 
for every odd $m$ and $P_m(1) = 2^{m/2}$ for every even $m$.

Our next lemma is stated as Lemma 3.5 in [E-16], where its proof may also be found.
It plays a key role in [E-19a] [E-19b], and [E-19c] as well. 

\proclaim{Lemma 2.2}
If $P_m$ and $Q_m$ are the $m$-th Rudin-Shapiro polynomials of degree $M-1$ with $M := 2^m$,
$\delta := \sin^2(\pi/8)$, and
$$z_j := e^{it_j}\,, \quad t_j := \frac{2\pi j}{M}\,, \qquad j \in {\Bbb Z}\,,$$
then
$$\max \{|P_m(z_j)|^2,|P_m(z_{j+1})|^2\} \geq \delta 2^{m+1} = 2\delta M\,.$$
\endproclaim

\proclaim{Lemma 2.3}
Using the notation of Lemma 2.2 we have
$$|P_m(e^{it})|^2 \geq \delta M\,, \qquad t \in \left[ t_j - \frac{\delta}{2M}, t_j + \frac{\delta}{2M} \right]\,,$$
for every $j \in {\Bbb Z}$ such that
$$|P_m(z_j)|^2 \geq \delta 2^{m+1} = 2\delta M\,.$$
\endproclaim

\demo{Proof}
The proof is a simple combination of the Mean Value Theorem and Bernstein's inequality (Lemma 3.4) 
applied to the (real) trigonometric polynomial of degree $M-1$ defined by $S(t) := P_m(e^{it})P_m(e^{-it})$. 
Recall that (2.1) implies $0 \leq S(t) = |P_m(e^{it})|^2 \leq 2M$ for every $t \in {\Bbb R}$.   
\qed \enddemo

Let, as before $M := 2^m$ with an odd $m$. We define
$$T(t) := \text {Re}((1 + e^{iMt} + e^{2iMt} + \cdots + e^{8iMt}) P_m(e^{it})) = 
\text {Re} \left( \frac{e^{9iMt}-1}{e^{iMt}-1} \, P_m(e^{it})\right) \,. \tag 2.2$$
Observe that $T$ is a real trigonometric polynomial of degree at most $\mu-1 := 9M-1$. For every sufficiently large natural 
number $n$ there is an odd integer $m$ such that 
$$2^{-75} \leq \gamma := \frac{\mu}{2n} = \frac{9 \cdot 2^m}{2n} < 2^{-72}\,. \tag 2.3$$ 
Observe that
$$\|T\| := \max_{t \in {\Bbb R}}{|T(t)|} = |T(0)| = 9|P_m(1)| = 9 \cdot 2^{(m+1)/2} = 9(2M)^{1/2}
= 6 \sqrt{\gamma n}\,. \tag 2.4$$ 

\proclaim{Lemma 2.4}
In the notation of Lemmas 2.2 and 2.3, for every $j \in {\Bbb Z}$ satisfying  
$$|P_m(z_j)|^2 \geq \delta 2^{m+1} = 2\delta M$$
there are
$$a_j \in \left[ t_j - \frac{3\pi}{32M},t_j - \frac{\pi}{32M}\right] \qquad \text{and} \qquad 
b_j \in \left[ t_j + \frac{\pi}{32M}, t_j + \frac{3\pi}{32M} \right]$$
such that
$$|T(a_j)| \geq (0.005) \|T\| = (0.03)\sqrt{\gamma n} \qquad \text{and} \qquad 
|T(b_j)| \geq (0.005) \|T\| = (0.03)\sqrt{\gamma n}\,.$$
\endproclaim

\demo{Proof}
We prove the statement about the existence of $b_j$ as the proof of the statement about the existence of $a_j$
is essentially the same. Let
$$P_m(e^{it}) = R(t) e^{i\alpha(t)}\,, \qquad R(t) = |P_m(e^{it})|\,,$$
where the function $\alpha$ could be chosen so that it is differentiable on any interval where
$P_m(e^{it})$ does not vanish. Then
$$ie^{it}P_m^\prime(e^{it}) = R^\prime(t) e^{i\alpha(t)} + R(t) e^{i\alpha(t)}(i\alpha^\prime(t))\,,$$
hence
$$\alpha^\prime(t) = \text{Re} \left( \frac{e^{it}P_m^\prime(e^{it})}{P_m(e^{it})} \right)$$
on any interval where $P_m(e^{it})$ does not vanish. Combining Bernstein's inequality (Lemma 3.4), Lemma 2.3, and
$\|P_m\| \leq (2M)^{1/2}$, we obtain
$$|\alpha^\prime(t)| \leq \frac{M(2M)^{1/2}}{(\delta M)^{1/2}} = 
\left( \frac{2}{\delta} \right)^{1/2} M \leq (3.7)M\,, \qquad t \in \left[ t_j, t_j + \frac{\delta}{2M} \right]\,. \tag 2.5$$
Now let
$$\frac{e^{9iMt}-1}{e^{iMt}-1} = \left| \frac{e^{9iMt}-1}{e^{iMt}-1} \right| e^{4Mt} \,, 
\qquad t \in \left( t_j -\frac{2\pi}{9M}, \, t_j + \frac{2\pi}{9M} \right)\,. \tag 2.6$$  
By writing 
$$(1 + e^{iMt} + e^{2iMt} + \cdots + e^{8iMt}) P_m(e^{it}) = \left| \frac{e^{9iMt}-1}{e^{iMt}-1} P_m(e^{it}) \right| e^{i(\alpha(t) + 4Mt)}\,,$$
we see by (2.5) and (2.6) that $\beta(t) := \alpha(t) + 4Mt$ satisfies
$$(0.3)M = 4M - (3.7)M \leq 4M - |\alpha^\prime(t)| \leq |\beta^\prime(t)|\,, 
\qquad t \in \left[t_j,t_j + \frac{\delta}{M} \right]\,. \tag 2.7$$
It is also simple to see that
$$\left| \frac{e^{9iMt}-1}{e^{iMt}-1} \right| \geq \, \left| \frac{e^{iM\pi}-1}{e^{iM\pi/9}-1} \right| = \frac{2}{2\sin(\pi/18)} 
\geq \frac{18}{\pi}\,, \qquad t \in \left[ t_j - \frac{\pi}{9M},t_j + \frac{\pi}{9M} \right]\,. \tag 2.8$$
Observe that (2.7) and (2.8) imply that there are
$$b_j \in \left[ t_j + \frac{\pi}{32M},t_j + \frac{3\pi}{32M} \right]$$
for which
$$\left| \frac{e^{9iMb_j}-1}{e^{iMb_j}-1} \right| \geq \frac{18}{\pi} \tag 2.9$$
and
$$\cos (\beta(b_j)) \geq \cos \left( \frac{\pi}{2} - \frac{(0.15)\pi}{16} \right) \geq 0.0294\,. \tag 2.10$$
Combining (2.9), (2.10), Lemma 2.3, and (2.4) we obtain
$$\split |T(b_j)| & = \, \left| \text {Re} \left( \frac{e^{9iMb_j}-1}{e^{iMb_j}-1} \, P_m(e^{ib_j})\right) \right|   
= \left| \frac{e^{9iMb_j}-1}{e^{iMb_j}-1} \right| \, \left| P_m(e^{ib_j}) \right| |\cos(\beta(b_j)| \cr 
& \geq \, \frac{18}{\pi} (\delta M)^{1/2} (0.0294) \geq \frac{(05292) \sin(\pi/8)}{9\sqrt{2}\pi} 9(2M)^{1/2}  
\geq (0.005) 9(2M)^{1/2} \cr 
& \geq \, (0.005)\|T\|\,. \endsplit $$
\qed \enddemo

\head 3. Tools from Approximation Theory \endhead

Let ${\Cal T}_\nu$ denote the set of all real trigonometric polynomials of degree at most $\nu$.
Let $\|T\|$ denote the maximum modulus of a trigonometric polynomial $T$ on ${\Bbb R}$.

\medskip

\noindent {\bf Definition 3.1} Let $n > 0$ be an integer divisible by $10$. We call ${\Cal I}$ {\it suitable} if

\medskip

\noindent (a) The endpoints of each interval in ${\Cal I}$ are in $(10\pi/n){\Bbb Z}$;

\noindent (b) ${\Cal I}$ is invariant under the maps $\theta \rightarrow \pi \pm \theta$;

\noindent (c) $|{\Cal I}| = 4N$ for some $N \leq \gamma n$.

\medskip

\noindent We call a suitable collection ${\Cal I}$ of disjoint intervals in ${\Bbb R}/(2\pi {\Bbb R})$ is
{\it well-separated} if

\medskip

\noindent (d) $|I| \leq 3990\pi/n$ for each $I \in {\Cal I}$;

\noindent (e) $d(I,J) \geq 10\pi/n$ for each $I,J \in {\Cal I}$ with $I \neq J$,

\noindent (f) The sets $\bigcup_{I \in {\Cal I}}{I}\,$ and $\,(\pi/2){\Bbb Z} + [-5\pi/n, 5\pi/n]$ are disjoint.

\medskip

We will denote the intervals in a suitable and well-separated collection ${\Cal I}$ by 
$$I_j, \qquad j=1,2,\ldots,4N\,,$$ 
where $I_1, I_2, \ldots, I_N \subset (0,\pi/2)$. 
Associated with an interval $[a,b] \subset [-\pi+5\pi/n,\pi-5\pi/n]$ we define
$$\Phi_{[a,b]}(t) :=
\cases
1, \quad & \text{if \enskip} t \in [a,b]\,,
\\
0, \quad & \text{if \enskip} t \in [-\pi,a-5\pi/n] \cup [b+5\pi/n,\pi]\,,
\\
(n/(5\pi))(t-a-5\pi/n), \quad &  \text{if \enskip} t \in [a-5\pi/n,a]\,,
\\
(n/(5\pi))((b+5\pi/n)-t), \quad & \text{if \enskip} t \in [b,b+5\pi/n]\,.
\endcases$$
We call the coloring $\alpha: {\Cal I} \rightarrow \{-1,1\}$ {\it symmetric} if 
$\alpha(I) = \alpha(\pi-I)$ and $\alpha(I) = -\alpha(\pi+I)$.  
Associated with a symmetric ${\Cal I}: \rightarrow \{-1,1\}$ let
$$g_{\alpha} := \sum_{j=1}^{4N}{\alpha(I_j) \Phi_{I_j}} \qquad \text{and} \qquad G_{\alpha} := K\sqrt{n}\,g_{\alpha}\,.$$
Let $S_o := \{1,3,\ldots,2n-1\}$ be the set of odd numbers between $1$ and $2n-1$. 
Let $C_{2\pi}$ denote the set of all continuous $2\pi$ periodic functions on ${\Bbb R}$.
Associated with $f \in C_{2\pi}$ we define the $n$th partial sum 
$$S_n(f,t) := a_0 + \sum_{k=1}^n{(a_k\cos(kt) + b_k \sin(kt))}$$
of the Fourier series expansion of $f$, where
$$a_0 = a_0(f) := \frac{1}{2\pi} \, \int_{-\pi}^{\pi}{f(t) \, dt}\,,$$
$$a_k = a_k(f) := \frac{1}{\pi} \, \int_{-\pi}^{\pi}{f(t) \cos(kt)\, dt}\,, \qquad k=1,2,\ldots\,,$$
and
$$b_k = b_k(f) := \frac{1}{\pi} \, \int_{-\pi}^{\pi}{f(t) \sin(kt)\, dt}\,, \qquad k=1,2,\ldots\,.$$
Observe that if $\alpha: {\Cal I} \rightarrow \{-1,1\}$ is symmetric, then
$$S_{2n}(G_\alpha,t) = S_{2n-1}(G_\alpha,t) = \sum_{k=1}^n{b_{2k-1}(G_\alpha) \sin((2k-1)t)}\,.$$
Associated with $f \in C_{2\pi}$ we also define
$$E_n(f) := \min_{Q \in {\Cal T}_n}{\|f-Q\|}$$
and
$$\omega(f,\delta) := \max_{t \in {\Bbb R}}{|f(t+\delta) - f(t)|}\,.$$

In the proof of Theorem 6.1 we will use D. Jackson's theorem on best uniform approximation
of continuous periodic functions with exact constant. The result below is due to Korneichuk [K-62].

\proclaim{Lemma 3.2} If $f \in C_{2\pi}$ then
$$E_n(f) \leq \omega\left( f,\frac{\pi}{n+1} \right)\,.$$
\endproclaim

In the proof of Theorem 6.1 we will also use the following result of De La Vall\'ee Poussin, 
the proof of which may be found on pages 273--274 in [D-93].  

\proclaim{Lemma 3.3} Associated with $f \in C_{2\pi}$ let
$$V_n(f,t) := \frac 1n \sum_{j=n}^{2n-1}{S_j(f,t)}\,.$$
We have
$$\max_{t \in {\Bbb R}}{|V_n(f,t) - f(t)|} \leq 4E_n(f)\,.$$
\endproclaim

The following inequality is known as Bernstein's inequality and plays an important role in the 
proof of Lemma 3.5.

\proclaim{Lemma 3.4} We have 
$$\|U^{(k)}\| \leq \nu^k \|U\|\,, \qquad U \in {\Cal T}_{\nu}\,, \qquad \nu = 1,2,\ldots\,, \quad k = 1,2,\ldots\,.$$
\endproclaim

\proclaim{Lemma 3.5}
Suppose $U \in {\Cal T}_\nu$, $\tau \in [0,2\pi/\nu]$, $A \geq 0.005$, and $|U(\tau)| \geq A \|U\|$.
Let
$$I_{j,\nu} := \left[ \frac{j\eta}{\nu}, \frac{(j+1)\eta}{\nu} \right] 
\subset \left[ \tau, \tau + \frac{18\pi}{\nu} \right]\,, \qquad j = u,u+1,\ldots, k\,. \tag 3.1$$
We have
$$\min_{t \in I_{j,\nu}}{|U(t)|} \geq \frac{A}{400} \left( \frac{\eta}{18\pi} \right)^{200} \|U\|$$
for at least one $j \in \{v,v+1,\ldots,v+399\}$ for every $v \in \{u,u+1,\ldots,k-399\}\,.$
\endproclaim

\demo{Proof}
Suppose the statement of the lemma is false, and there are $v \in \{u,u+1,\ldots,k-399\}$ and
$$x_j \in I_{j,\nu} := \left[ \frac{j\eta}{\nu}, \frac{(j+1)\eta}{\nu} \right] 
\subset \left[ \tau, \tau + \frac{18\pi}{\nu} \right] \tag 3.2$$
such that
$$|U(x_j)| < \frac{A}{400} \left( \frac{\eta}{2\pi} \right)^{200} \|U\|\,, \qquad j \in \{v,v+1,\ldots,v+399\}\,.$$
Let $y_j := x_{v+2j-1}$ for $j \in \{1,2,\ldots,200\}$. Then the points $y_j$ satisfy
$$y_1 - \tau \geq \frac{\eta}{\nu}\, \quad \text{and} \quad y_{j+1}-y_j \geq \frac{\eta}{\nu}\,, 
\qquad j \in \{1,2,\ldots 200\}\,.$$
By the well-known formula for divided differences we have
$$U(\tau) \prod_{h=1}^{200}{(\tau-y_h)^{-1}} + 
\sum_{j=1}^{200}{U(y_j) (\tau-y_j)^{-1} \prod_{h=1 \atop h \neq j}^{200}{(y_h-y_j)^{-1}}} = 
\frac{1}{200!} U^{(200)}(\xi)\,,$$
and combining this with $|U(\tau)| \geq A \|U\|$, (3.1), and (3.2), we get 
$$ A \|U\| \left( \frac{18\pi}{\nu} \right)^{-200} \leq 
200 \, \frac{A}{400} \left( \frac{\eta}{18\pi} \right)^{200} \|U\| \left( \frac {\eta}{\nu} \right)^{-200} 
+ \frac{1}{200!} |U^{(200)}(\xi)|\,,$$
with some $\xi \in [\tau, \tau + 2\pi/\nu]$. Therefore Bernstein's inequality (Lemma 3.4) yields that
$$A \|U\| \left( \frac{18\pi}{\nu} \right)^{-200} 
\leq 200 \, \frac{A}{400} \left( \frac{\eta}{18\pi} \right)^{200} \|U\| \left( \frac {\eta}{\nu} \right)^{-200}
+ \frac{1}{200!} \nu^{200} \|U\|\,,$$
that is,
$$A \leq \frac{2(18\pi)^{200}}{200!} \leq 2 \left( \frac{18\pi e}{200} \right)^{200} < 0.005\,,$$
which contradicts our assumption $A \geq 0.005$.
\qed \enddemo

The following lemma ascribed to M. Riesz is well-known and can easily be proved by a simple
zero counting argument (see [B-95], for instance).

\proclaim{Lemma 3.6}
If $T \in {\Cal T}_\nu$, $t_0 \in {\Bbb R}$, and $|T(t_0)| = \|T(t)\|$, then
$$|T(t)| \geq |T(t_0)| \, \cos(\nu(t-t_0))\,, \qquad t \in {\Bbb R}, \enskip |t-t_0| \leq \frac{\pi}{2\nu}\,.$$
\endproclaim

We will also need the following simple corollary of the above lemma.

\proclaim{Lemma 3.7}
If $L  = 32n$,
$$t_r := \frac{(2r-1)\pi}{4L}\,, \qquad r=1,2,\ldots,4L\,,$$
and $T \in {\Cal T}_n$, then
$$\max_{t \in {\Bbb R}}{|T(t)|} \leq (\cos(\pi/64))^{-1} \max_{1 \leq r \leq 4L}{|T(t_r)|} 
\leq (1.0013) \, \max_{1 \leq r \leq 4L}{|T(t_r)|}\,.$$
\endproclaim

\head 4. Minimizing Discrepancy \endhead

Associated with a vector $\bold{x} = \langle x_1,x_2,\ldots,x_v \rangle \in {\Bbb R}^v$ let 
$$\|\bold{x}\|_\infty := \max\{|x_1|,|x_1|,\ldots,|x_v|\}\,.$$ 
A crucial ingredient in [B-20] is the main ``partial coloring" lemma of Spencer [S-85]
based on a technique of Beck [B-81]. In Section 4 of [B-20] a simple consequence of a
variant of this due to Lovett and Meka [L-15, Theorem 4] is observed, and it plays
an important part in the proof of Theorem 6.1. This can be stated as follows.

\proclaim{Lemma 4.1}
Let $\bold{y}_1, \bold{y}_2,\ldots,\bold{y}_u \in {\Bbb R}^v$ and $\bold{x}_0 \in [-1,1]^v$. If 
$c_1,c_2,\ldots,c_u \geq 0$ are such that 
$$\sum_{r=1}^u\exp(-(c_r/14)^2) \leq \frac{v}{16}\,, \tag 4.1$$
then there exists an $\bold{x} \in \{-1,1\}^v$ such that 
$$|\langle \bold{x} - \bold{x}_0, \bold{y}_r \rangle| \leq (c_r + 30)\sqrt{u} \, \|\bold{y}_r\|_\infty\,, 
\qquad r=1,2,\ldots,u\,.$$
\endproclaim

\head 5. The Cosine Polynomial \endhead

\proclaim{Theorem 5.1}
Let $n > 0$ be a sufficiently large integer divisible by $10$. There exist a cosine polynomial 
$$c(t) = \sum_{k=0}^{\mu}{\varepsilon_k\cos(2kt)}\,, \qquad \varepsilon_k \in \{-1,1\}\,, 
\quad k=1,2,\ldots,\mu\,, \tag 5.1$$
and a suitable and well-separated collection ${\Cal I}$ of disjoint intervals in ${\Bbb R}/(2\pi {\Bbb Z})$ such that 
$$c(t) \geq \eta_1 \sqrt{n}\,, \qquad t \notin \bigcup_{I \in {\Cal I}}{I}\,,$$
and
$$c(t) \leq \sqrt{n}\,, \qquad t \in {\Bbb R}\,,$$ 
where $\eta_1 > 0$ is an absolute constant. 
\endproclaim

\demo{Proof}
Let $c(t) := U(t) := T(2t)$, where $T \in {\Cal T}_\mu$ with $\mu := 9M$ is defined by (2.2) and 
$U \in {\Cal T}_{\nu}$ with $\nu := 2\mu$. Observe that 
$c$ is of the form (5.1). It follows from (2.1), (2.3), and $2^{-75} < \gamma \leq 2^{-72}$ that  
$$|c(t)| \leq 9\sqrt{2M} \leq 3\sqrt{2\mu} \leq \sqrt{n}\,.$$
Set 
$$\eta := 10\pi \gamma = 10\pi(2\mu/n) \qquad \text{and} \qquad 
\eta_1 := \frac{0.005}{400} \left( \frac{\eta}{18\pi} \right)^{200}\,.$$ 
We partition ${\Bbb R}/(2\pi {\Bbb Z})$ into $n/5$ intervals 
$$I_j := [10\pi j/n,10\pi(j+1)/n]\,, \qquad j=0,1,\ldots,n/5-1\,,$$ 
and say that an interval $I_j$ is good if 
$$\min_{t \in I_j}{|U(t)|} \geq \frac{0.005}{400} \left( \frac{\eta}{18\pi} \right)^{200} \|U\|\,.$$
Let ${\Cal J}$ be the collection of maximal unions of consecutive good intervals $I_j$, and 
let ${\Cal I}$ be the collection of the remaining intervals (that is, the maximal unions of 
consecutive bad intervals). We claim that ${\Cal I}$ is the required suitable and well-separated 
collection. 

First, to see that ${\Cal I}$ is suitable, note that the endpoints of each of the intervals $I_j$ are 
in $10\pi{\Bbb Z}$. The set ${\Cal I}$ is invariant under the maps $\theta \rightarrow \pi \pm \theta$ 
by the symmetries of the functions $\cos(2kt)$, $k=0,1,\ldots,\mu$. To see that $4N = |{\Cal I}| \leq 4\gamma n$, 
note that a real trigonometric polynomial of degree at most $\nu$ has at most $2\nu$ real zeros 
in a period, and hence there are at most $4\nu$ values of $t$ in a period for which 
$$U(t) = \frac{\pm 0.005}{400} \left( \frac{\eta}{18\pi} \right)^{200} \|U\|\,.$$
Since each $I \in {\Cal I}$ must contain at least two such points (counted with multiplicities), 
we have $4N := |{\Cal I}| \leq 2\nu = 4\gamma n$. Thus ${\Cal I}$ has each of the properties (a), 
(b) and (c) in the definition of a suitable collection. 

We now show that ${\Cal I}$ is well-separated. By Lemmas 3.5 and 2.4 any $400$ consecutive intervals 
$I_j$ must contain a good interval, and hence $|I| \leq 3990\pi/n$ for each $I \in {\Cal I}$. 
Thus ${\Cal I}$ has property (d) in the definition of a well-separated collection. 
The fact that ${\Cal I}$ has property (e) in the definition of a suitable collection 
is obvious by the construction. Finally observe that for an even $m$ we have 
$$|P_m(1)| = 2^{(m+1)/2} = \|P_m(e^{it})\|\,,$$ 
from which 
$$|T(0)| = |T(\pi)| = \|T\|$$
follows. Hence, property (f) in the definition of a well-separated collection 
follows from the Riesz's Lemma stated as Lemma 3.6 (recall that $\nu = 2\mu = \gamma n < 2^{-72}n$).     
\qed \enddemo

\head 6. The Sine Polynomials \endhead

\proclaim{Theorem 6.1}
Let $n > 0$ be an integer divisible by $10$.
Let ${\Cal I}$ be a suitable and well-separated collection of disjoint intervals in ${\Bbb R}/(2\pi {\Bbb Z})$. 
There exists a sine polynomial
$$s_o(t) = \sum_{k=1}^n{\varepsilon(2k-1)\sin(2k-1)t)}\,, \qquad \varepsilon(2k-1) \in \{-1,1\}\,,$$
such that
$$|s_o(t)| \geq 36 \sqrt{n}\,, \qquad t \in \bigcup_{I \in {\Cal I}}{I}\,, \qquad \text{and} \qquad
|s_o(t)| \leq 1090\sqrt{n}\,, \qquad t \in {\Bbb R}\,.$$ 
\endproclaim

To prove Theorem 6.1 we need some lemmas.

\proclaim{Lemma 6.2} Let ${\Cal I}$ be a suitable and well-separated collection of disjoint intervals in 
${\Bbb R}/(2\pi {\Bbb R})$. There exists a symmetric coloring $\alpha: {\Cal I} \rightarrow \{-1,1\}$ 
such that 
$$a_k(G_\alpha) = 0\,, \qquad k=0,1,\ldots,2n\,,$$
$$b_{2k}(G_\alpha) = 0\,, \qquad \text{and} \qquad |b_{2k-1}(G_\alpha)| \leq 1 \,, \qquad k=1,2,\ldots,n\,.$$ 
\endproclaim

\demo{Proof}
As before, we denote the intervals in a suitable and well-separated collection ${\Cal I}$ by
$I_j$, $j=1,2,\ldots,4N$, where $I_1, I_2, \ldots, I_N \subset (0,\pi/2)$. As we have already observed before, we have 
$a_k(G_\alpha) = 0, \enskip k=0,1,\ldots,2n$, and $b_{2k}(G_\alpha) = 0, \enskip k=1,2,\ldots,n$, for every symmetric coloring 
$\alpha: {\Cal I} \rightarrow \{-1,1\}$, so we have to show only that there exists a symmetric coloring $\alpha: {\Cal I} \rightarrow \{-1,1\}$ 
such that $|b_{2k-1}(G_\alpha)| \leq 1, \enskip k=1,2,\ldots,n$. To this end let 
$$\bold{y}_k := \langle y_{k,1},y_{k,2},\ldots,y_{k,N} \rangle\,, \qquad k=1,2,\ldots,n\,,$$
with
$$y_{k,j} := \frac{4K\sqrt{n}}{\pi} \int_{-\pi}^{\pi} {\Phi_{I_j}(t)\sin((2k-1)t) \, dt}\,, \qquad k=1,2,\ldots,n, \enskip j=1,2,\ldots,N\,.$$ 
If $\alpha: {\Cal I} \rightarrow \{-1,1\}$ is a symmetric coloring, then by the symmetry conditions on ${\Cal I}$ we have
$$b_{2k-1}(G_\alpha) := \frac{1}{\pi} \int_{-\pi}^{\pi}{G_\alpha(t) \sin((2k-1)t)) \, dt}
= \sum_{j=1}^N{\alpha(I_j)y_{k,j}}\,, \qquad k=1,2,\ldots,n\,.$$
We apply Lemma 4.1 with $u := n$, $v := N$, $\bold{x}_0 := \bold{0} \in [-1,1]^N$, and                                      
$$c_1 = c_2 = \cdots = c_n := 14\sqrt{\log(16n/N)}\,.$$
Observe that
$$\sum_{r=1}^u{\exp(-c_r^2/14^2)} = n \frac{N}{16n} = \frac{N}{16}\,,$$
so (4.1) is satisfied. It follows from Lemma 4.1 that there exists an 
$$\langle \alpha(I_1),\alpha(I_2),\ldots,\alpha(I_N) \rangle = \bold{x} \in \{-1,1\}^N$$ 
such that
$$|\langle \bold{x},\bold{y}_k \rangle| \leq (c_k+30)\sqrt{N} \, \|\bold{y}_k\|_\infty\,, \qquad k=1,2,\ldots,n\,.$$ 
As ${\Cal I}$ is well-separated, by part (d) of the definition we have
$$|y_{k,j}| \leq \frac{4K \sqrt{n}}{\pi}(|I_j| + 10/\pi/n) \leq \frac{4K \sqrt{n}}{\pi}\frac{4000\pi}{n} = \frac{16000K}{\sqrt{n}}$$ 
for every $k = 1,2,\ldots,n$ and $j = 1,2,\ldots,N\,.$
It follows that
$$|b_{2k-1}(G_\alpha)| = |\langle \bold{x},\bold{y}_k \rangle| \leq (14\sqrt{\log(16n/N)} + 30) \sqrt{N/n} \cdot 16000 K\,, 
\qquad k=1,2,\ldots,n\,.$$
As the right-hand side above is an increasing function of $N$ for $N/n \leq \gamma < 1$, we have 
$$|b_{2k-1}(G_\alpha)| = |\langle \bold{x},\bold{y}_k \rangle| \leq (14\sqrt{\log(16/\gamma} + 30) \sqrt{\gamma} \cdot 16000 K \leq 1\,,
\qquad k=1,2,\ldots,n\,,$$
where the last inequality follows from $K := 2^9$ and the inequality $2^{-75} \leq \gamma < 2^{-72}$. 
Hence the desired symmetric coloring is given by setting 
$$\langle \alpha(I_1),\alpha(I_2),\ldots,\alpha(I_N) \rangle := \bold{x}\,.$$  
\qed \enddemo

From now on let $\alpha:{\Cal I} \rightarrow \{-1,1\}$ denote the symmetric coloring 
guaranteed by Lemma 6.2. We have
$$V_n(G_\alpha,t) = \sum_{k=1}^n{\widetilde{\varepsilon}(2k-1)\sin((2k-1)t)}\,, \qquad 
|\widetilde{\varepsilon}(2k-1)| \leq 1\,.$$

\proclaim{Lemma 6.3} There is a coloring $\varepsilon: S_o \rightarrow \{-1,1\}$ such that with the notation 
$$s_o(t) = \sum_{k=1}^n{\varepsilon(2k-1)\sin((2k-1)t)}$$
we have
$$|s_o(t) - V_n(G_\alpha,t)| \leq 66\sqrt{n}\,, \qquad t\in {\Bbb R}\,.$$
\endproclaim

\demo{Proof}
Let $L := 32n$, 
$$t_r := \frac{(2r-1)\pi}{4L}\,, \qquad r=1,2,\ldots,4L\,,$$
$$y_{r,k} := \sin((2k-1)t_r)\,, \qquad r=1,2,\ldots,L\,, \enskip k=1,2,\ldots,n\,,$$
$$\bold{y}_r := \langle y_{r,1}, y_{r,2}, \ldots,y_{r,n} \rangle\,, \qquad r=1,2,\ldots,L\,.$$
Observe that 
$$s_o(t_r) - V_n(G_\alpha,t_r) = \sum_{k=1}^n{(\varepsilon(2k-1) - \widetilde{\varepsilon}(2k-1))y_{r,k}}
= \langle \bold{e} - \widetilde{\bold{e}}, \bold{y}_r \rangle\,, \tag 6.1$$
where
$$\bold{e} := \langle \varepsilon(1), \varepsilon(3),\ldots, \varepsilon(2n-1) \rangle \qquad \text{and} \qquad 
\widetilde{\bold{e}} := 
\langle \widetilde{\varepsilon}(1), \widetilde{\varepsilon}(3),\ldots, \widetilde{\varepsilon}(2n-1) \rangle\,.$$
We apply Lemma 4.1 with $u := L$, $v := n$, $\bold{x}_0 := \widetilde{\bold{e}}$, and 
$$c_1 = c_2 = \cdots = c_n := 42\sqrt{\log 2}\,.$$ 
Observe that
$$\sum_{r=1}^u{\exp(-c_r^2/14^2)} = L2^{-9} = \frac{n}{16}\,,$$
so (4.1) is satisfied. It follows from Lemma 4.1 that there exists an $\bold{e} \in \{-1,1\}^n$ 
such that
$$|\langle \bold{e} - \widetilde{\bold{e}}, \bold{y}_r \rangle| \leq (c_r + 30)\sqrt{n} \|\bold{y}_r\|_\infty 
\leq (c_r + 30)\sqrt{n} \leq 65\sqrt{n}\,,\,, \qquad r=1,2,\ldots,L\,. \tag 6.2$$
Combining (6.1) and (6.2) we obtain
$$|s_o(t_r) - V_n(G_\alpha,t_r)| \leq 65 \sqrt{n}\,, \qquad r=1,2,\ldots,L\,.$$ 
Note that by the special form of the trigonometric polynomials $s_o$ and $V_n(G_\alpha,\cdot)$ we have 
$$\max_{1 \leq r \leq L}{|s_o(t_r) - V_n(G_\alpha,t_r)|} = \max_{1 \leq r \leq 4L}{|s_o(t_r) - V_n(G_\alpha,t_r)|}\,,$$
hence
$$|s_o(t_r) - V_n(G_\alpha,t_r)| \leq 65 \sqrt{n}\,, \qquad r=1,2,\ldots,4L\,.$$
This, together with Lemma 3.7 gives the lemma. 
\qed \enddemo

\proclaim{Lemma 6.4}
We have
$$|V_n(G_\alpha,t)| \geq \frac{K\sqrt{n}}{5}\,, \qquad t \in \bigcup_{I \in {\Cal I}}{I}\,, \qquad \text{and} \qquad 
|V_n(G_\alpha,t)| \leq 2K\sqrt{n}\,, \qquad t \in {\Bbb R}\,.$$
\endproclaim

\demo{Proof}
Combining Lemma 3.3 and 3.2 we have 
$$\max_{t \in {\Bbb R}}{|V_n(G_\alpha,t) - G_\alpha(t)|} \leq 4E_n(G_\alpha)
\leq 4\omega(G_\alpha,\pi/n) \leq \frac{4K\sqrt{n}}{5}\,,$$
and the lemma follows.
\qed \enddemo

Let 
$$s_e(t) := \text{Im}(P_{<(n+1)}(e^{2it})) - \text{Im}(P_{<(\nu+1)}(e^{2it}))\,.$$

\proclaim{Lemma 6.5}
We have 
$$\|s_e\| \leq 6\sqrt{n}\,.$$
\endproclaim

\demo{Proof}
This is an obvious consequence of Lemma 2.1. Recall that $\nu = \gamma n \leq 2^{-72}n$.
\qed \enddemo

\demo{Proof of Theorem 6.1}
Let ${\Cal I}$ be a suitable and well-separated collection of disjoint intervals in ${\Bbb R}/(2\pi {\Bbb Z})$. By 
Lemma 6.3 there is a coloring $\varepsilon:S_o \rightarrow \{-1,1\}$ such that if $\alpha: {\Cal I} \rightarrow \{-1,1\}$
is the symmetric coloring given by Lemma 6.2, then
$$|s_o(t) - V_n(G_\alpha,t)| \leq 66\sqrt{n}\,, \qquad t\in {\Bbb R}\,.$$ 
Hence by Lemma 6.4 and $K := 2^9$ we have
$$|s_o(t)| \geq |V_n(G_\alpha,t)| - |s_o(t) - V_n(G_\alpha,t)|  
\geq 102\sqrt{n} - 66\sqrt{n} \geq 36\sqrt{n}\,, \qquad t \in \bigcup_{I \in {\Cal I}}{I}\,,$$
and
$$|s_o(t)| \leq |V_n(G_\alpha,t)| + |s_o(t) - V_n(G_\alpha,t)|  
\leq 2^{10}\sqrt{n} + 66\sqrt{n} \leq 1090\sqrt{n}\,, \quad t \in {\Bbb R}\,.$$
\qed \enddemo

\head 7. Proof of Theorems 1.1 and 1.2 \endhead

\demo{Proof of the Theorems 1.2} 
It is sufficient to prove the theorem with $2n$ replaced by $4n$ and without loss of generality we may assume 
that $n > 0$ is an integer divisible by $10$. Since the Littlewood polynomial $P_{4n}(z) := 1-z-z^2-\cdots-z^{4n}$ 
does not vanish on the unit circle, we may assume also that $n$ is sufficiently large.
By Theorems 5.1 and 6.1 the Littlewood polynomial $P_{4n}$ of degree $4n$ defined by
$$P_{4n}(e^{it})e^{-2int} = (-1 + 2c(t)) + 2i(s_o(t) + s_e(t))$$
has the properties required by the theorem. It is obvious from the construction that the coefficients of $P_{4n}$
satisfy the requirements. To see that the required inequalities are satisfied let ${\Cal I}$ be a suitable and well-separated 
collection of disjoint intervals in ${\Bbb R}/(2\pi {\Bbb Z})$ on which (5.1) holds. Then Theorem 5.1 gives that
$$|P_{4n}(e^{it})| \geq |-1 + 2c(t)| \geq \eta_1 \sqrt{n}\,, \qquad t \notin \bigcup_{I \in {\Cal I}}{I}\,,$$ 
while Theorem 6.1 gives that
$$|P_{4n}(e^{it})| \geq |2(s_o(t) + s_e(t))| \geq |2s_o(t)| - |2s_e(t)| \geq 72\sqrt{n} - 12\sqrt{n}  
= 60\sqrt{n}\,, \quad t \in \bigcup_{I \in {\Cal I}}{I}\,.$$
Combining the two inequalities above gives the lower bound of the theorem. The upper bounds of the theorem follows 
from combining the upper bounds of Theorems 5.1 and 6.1 by
$$\split |P_{4n}(e^{it})| & \leq \, |-1 + 2c(t)| + |2(s_o(t) + s_e(t))| \leq 1 + 2\sqrt{n} + 2180\sqrt{n} + 12 \sqrt{n} \cr
& \leq \, 1 + 2196 \sqrt{n}\,, \quad {t \in \Bbb R}\,. \cr \endsplit$$ 
For the value $m_n$ in the theorem we have $m_n = 2\mu =  2\gamma n$, so $\eta = 2\gamma > 0$ can be chosen. 
\qed \enddemo

\enddocument